\documentclass{amsart}
\usepackage{amssymb, amsmath, amsfonts, amscd}

\input xypic

\theoremstyle{plain}
\newtheorem{theorem}{Theorem}[section]

\newtheorem{proposition}[theorem]{Proposition}

\newtheorem{example}[theorem]{Example}

\theoremstyle{definition}
\newtheorem{definition}[theorem]{Definition}

\theoremstyle{remark}

\numberwithin{equation}{theorem}

\renewcommand{\O}{\mathcal{O} }

\renewcommand{\P}{\mathbf{P} }
\renewcommand{\Pr}{\mathcal{P} }

\newcommand{\vw}{\mathbf{V}_W}
\newcommand{\pa}{\partial}
\newcommand{\onab}{\overline{\nabla}}

\newcommand{\Der}{\operatorname{Der} }
\newcommand{\Hom}{\operatorname{Hom} }
\newcommand{\Ext}{\operatorname{Ext} }
\newcommand{\End}{\operatorname{End} }
\newcommand{\Spec}{\operatorname{Spec} }

\renewcommand{\H}{\operatorname{H} }

\newcommand{\D}{\operatorname{D} }
\newcommand{\HH}{\operatorname{HH} }

\newcommand{\C}{\operatorname{C} }

\newcommand{\U}{\operatorname{U} }
\renewcommand{\lg}{\mathfrak{g}}
\begin{document}

\title{Chern classes and Lie-Rinehart algebras}
\author{Helge Maakestad }
\address{Emmy Noether Institute for Mathematics, Israel}
\email{makesth@macs.biu.ac.il }
\thanks{This work was partially supported by the Emmy Noether Research Institute for
Mathematics (center of the Minerva Foundation of Germany), the Excellency Center ``Group
Theoretic Methods in the Study of Algebraic Varieties''  
and EAGER (EU network HPRN-CT-2009-00099) }
\date{march 2002}
\begin{abstract} Let $A$ be a $F$-algebra where $F$ is a field, and let
$W$ be an $A$-module of finite presentation. We use the linear  Lie-Rinehart algebra
$\vw$ of $W$ to define the first Chern-class $c_1(W)$ in $\H^2(\vw|_U,\O_U)$,
where $U$ in $\Spec(A)$ is the open subset where $W$ is locally free. 
We compute explicitly algebraic $\vw$-connections on maximal 
Cohen-Macaulay modules $W$ on the hypersurface-singularities 
$B_{mn2}=x^m+y^n+z^2$,
and show that these connections are integrable, hence the first Chern-class
$c_1(W)$ vanishes. We also look at indecomposable maximal 
Cohen-Macaulay modules on  quotient-singularities in dimension 2, 
and prove that their first Chern-class vanish.
\end{abstract}
\maketitle

\tableofcontents

\section*{Introduction} 
Classically, the Chern-classes of a locally free coherent $A$-module $W$ are defined using the curvature
$R_{\nabla}$ of a connection $\nabla:W\rightarrow W\otimes \Omega^1_X$. A connection
$\nabla$ gives rise to a covariant derivation $\overline{\nabla}:\Der_F(A)\rightarrow
\End_F(W)$. If we more generally consider a coherent $A$-module $W$, a connection $\nabla$
might not exist. In this paper we will consider the problem of defining Chern-classes
for $A$-modules $W$ using a covariant derivation defined on a certain sub-Lie-algebra $\vw$
of $\Der_F(A)$: For an arbitrary $A$-module $W$ there exists a sub-Lie-algebra
$\vw$ of $\Der_F(A)$ called the linear Lie-Rinehart algebra, and also the notion of a $\vw$-connection.
%
%
There exists a complex,
the Chevalley-Hochshild complex $\C^\bullet(\vw,W)$ for  the $A$-module $W$ with
a flat $\vw$-connection, generalizing the classical deRham-complex.
If $A$ is a regular ring, the derivation module $\Der_F(A)$ is locally free, and
it follows that the complex $\C^\bullet(\Der_F(A),A)$ is quasi-isomorphic to the complex
$\Omega^\bullet_{A/F}$, hence the Chevalley-Hochschild complex of $\Der_F(A)$ can be 
used to compute the algebraic deRham-cohomology of $A$. The complex
$\C^\bullet(\vw,A)$ generalizes simultaneously the algebraic deRham-complex  and 
the Chevalley-Eilenberg complex.
This is due to \cite{PALAIS}. A natural thing to do is 
to investigate 
possibilities of defining Chern-classes of $A$-modules equipped with a $\vw$-connection,
generalizing the classical Chern-classes defined using the curvature of a connection.
Invariants for  Lie-Rinehart algebras have been considered
by several authors (see \cite{FER}, \cite{HUEBSCH} and \cite{KUB}), however
invariants for $A$-modules with a $\vw$-connection, where $\vw$ is the linear  Lie-Rinehart algebra
of $W$ does not appear to be treated in the litterature and that is the aim of this work. 
In this paper we also develop techniques to do explicit calculations of Chern-classss
of maximal Cohen-Macaulay modules on hypersurface singularities and two-dimensional
quotient singularities. 
 
We define for any $F$-algebra $A$  where $F$ is any
field, and any
$A$-module $W$ which is locally free on an open subset $U$ of $\Spec(A)$, the first
Chern-class $c_1(W)$ in $\H^2(\vw|_U,\O_U)$, where $\H^2(\vw|_U,\O_U)$ is the 
Chevalley-Hochschild cohomology of the restricted linear  Lie-Rinehart algebra $\vw$ of $W$, with 
values in the sheaf $\O_U$. This is Theorem \ref{chernclass}. We also 
prove in Theorem \ref{brieskorn} existence of explicit $\vw$-connections
$\overline{\nabla}^{\psi,\phi}$ on a class of maximal Cohen-Macaulay modules $W$ on the 
Brieskorn singularities $B_{mn2}$, which in fact are defined over any field $F$ of characteristic 
prime to $m$ and $n$.  We prove in Theorem \ref{zero} that
the $\vw$-connections defined in Theorem \ref{brieskorn} are all regular, hence the first
Chern-class is zero. Finally we prove in Theorem \ref{quotient}
that for any maximal Cohen-Macaulay module $W_\rho$ on any two-dimensional quotient-singularity
$\mathbf{C}^2/G$, the first chern class $c_1(W_\rho)$ is zero.

\section{Kodaira-Spencer maps and linear Lie-Rinehart algebras}

Let $A$ be an $F$-algebra, where $F$ is any field, and let $W$ be an $A$-module.
in this  section we use the Kodaira-Spencer class and the Kodaira-Spencer map to
define the linear Lie-Rinehart algebra $\vw$ of $W$, and the obstruction $lc(W)$ in 
$\Ext^1_A(\vw,\End_A(W))$ for existence of a $\vw$-connection on $W$. 

\begin{definition} Let $P$ be an $A$-bimodule. The \emph{Hochschild-complex}
of $A$ with values in $P$ $\C^\bullet(A,P)$ is defined as follows:
\[ \C^p(A,P)=\Hom_F(A^{\otimes p},P) \]
with differentials $d^p:\C^p(A,P)\rightarrow \C^{p+1}(A,P)$
defined by
\[ d^p\phi(a_1\otimes \cdots \otimes a_p\otimes a_{p+1})=
a_1\phi(a_2\otimes \cdots \otimes a_{p+1})+ \]
\[ \sum_{1\leq i\leq p}(-1)^i\phi(a_1\otimes \cdots \otimes a_ia_{i+1}
\otimes \cdots \otimes a_{p+1})+(-1)^{p+1}\phi(a_1\otimes \cdots 
\otimes a_p)a_{p+1} .\]
We adopt the convention that $\C^0(A,P)=P$, and $d^0(p)(a)=pa-ap$
for all $p$ in $P$, and $a$ in $A$. The i'th cohomology 
$\H^i(\C^\bullet(A,P))$ is denoted $\HH^i(A,P)$.
\end{definition}
There exists an exact sequence 

\begin{equation} \label{ex1}
0\rightarrow \HH^0(A,P)\rightarrow P \rightarrow \Der_F(A,P)
\rightarrow \HH^1(A,P)\rightarrow 0,
\end{equation}
and it is well known that if we let $P=\Hom_F(W,V)$ where $W$ and $V$ are 
$A$-modules, then
$\HH^i(A,P)$ equals $\Ext^i_A(W,V)$. Recall the definition of a
\emph{connection} on $W$: it is an $F$-linear map
$\nabla:W\rightarrow W\otimes \Omega^1_{A/F}$ with the property that
$\nabla(aw)=a\nabla(w)+w\otimes da$, where $d:A\rightarrow \Omega^1_{A/F}$
is the \emph{universal derivation}.
Put $P=\Hom_F(W,W\otimes \Omega^1_{A/F})$ in \ref{ex1}, and construct 
an element $C$ in $\Der_F(A,\Hom_F(W,W\otimes \Omega^1_{A/F}))$ in the 
following way: $C(a)(w)=w\otimes da$.

\begin{definition} The class $\overline{C}=ks(W)$ in 
$\Ext^1_A(W,W\otimes \Omega^1_{A/F})$ is the \emph{Kodaira-Spencer class}
of $W$.
\end{definition}
Note that $ks(W)$ is also referred to as the \emph{Atiyah-class} of $W$.

\begin{proposition} Let $A$ be any $F$-algebra, and let $W$ be any
$A$-module, then $ks(W)=0$ if and only if $W$ has a connection.
\end{proposition}
\begin{proof} We see that $ks(W)=0$ if and only if there exists an element
\[ \nabla: W \rightarrow W\otimes \Omega^1_{A/F}\]
 with the property that
 $d^0\nabla=C$. This is if and only if 
\[ (\nabla a-a\nabla)(w)=\nabla(aw)-a\nabla(w)=C(a)(w)=w\otimes da ,\]
hence $\nabla$ is a connection, and the claim follows.
\end{proof}
\begin{definition} \label{liealgebroid} Let $A$ be an $F$-algebra where $F$ is any
  field. A
\emph{ Lie-Rinehart algebra} on $A$ is a $F$-Lie-algebra and an $A$-module
$\mathfrak{g}$ with a map $\alpha:\lg\rightarrow \Der_F(A)$ satisfying
the following properties:
\begin{align} 
&\alpha(a\delta)=a\alpha(\delta)   \\
&\alpha([\delta, \eta])=[\alpha(\delta),\alpha(\eta)] \\
&[\delta,a\eta]=a[\delta,\eta]+\alpha(\delta)(a)\eta 
\end{align}
for all $a\in A$ and $\delta,\eta \in \lg$.
Let $W$ be an $A$-module. A $\lg$-\emph{connection} $\nabla$ on $W$,
is
an $A$-linear map $\nabla:\lg\rightarrow \End_F(W)$ wich satisfies
the \emph{Leibniz-property}, i.e. 
\[ \nabla (\delta)(aw)=a\nabla (\delta)(w)+\alpha (\delta)(a)w \]
for all $a\in A$ and $w\in W$. 
We say that $(W,\nabla) $ is a \emph{$\lg$-module}
if $\nabla$ is a homomorphism of Lie-algebras. The \emph{curvature} of the 
$\lg$-connection, $R_\nabla$ is defined as follows: 
\[ R_\nabla(\delta \wedge \eta)=[\nabla_\delta,\nabla_\eta]-\nabla_{[\delta,\eta]}.\]
\end{definition}

\begin{example}
Any connection $\nabla$ on $W$, gives an action $\nabla:\Der_F(A)\rightarrow
\End_F(W)$ with the property that $\nabla(\delta)(aw)=a\nabla(\delta)(w)+
\delta(a)w$ for any $\delta$ in $\Der_F(A)$, $a$ in $A$ and $w$ in $W$.
The Lie-algebra $\lg=\Der_F(A)$ is trivially a  Lie-Rinehart algebra, hence
$W$ has a $\lg$-connection. 
\end{example}

Given any Lie-algebroid $\lg$, and any $A$-module $W$ with a $\lg$-connection, the set of
$\lg$-connections on $W$ is a \emph{torsor} on the set $\Hom_A(\lg,\End_A(W))$.
Put $P=\Hom_F(W,W)$ in \ref{ex1}, and define for all $\delta$ in $\Der_F(A)$
the following element $C(\delta)$ in $\Der_F(A,\Hom_F(W,W))$:
$C(\delta)(a)(m)=\delta(a)m$. We get an $A$-linear map
\[ C:\Der_F(A)\rightarrow \Der_F(A,\Hom_F(W,W)).\]
\begin{definition} Let $A$ be any $F$-algebra, and let $W$ be any $A$-module.
We define the \emph{Kodaira-Spencer map}
\[ g:\Der_F(A)\rightarrow \Ext^1_A(W,W)\]
as follows: $g(\delta)=\overline{C(\delta)}$ in sequence \ref{ex1}.
We let $ker(g)=\vw$ be the \emph{linear  Lie-Rinehart algebra} of $W$.
\end{definition}
One immediately checks that the $A$-sub module $\vw$ of $\Der_F(A)$ satisfies
the axioms of definition \ref{liealgebroid}, hence $\vw$ is indeed
a  Lie-Rinehart algebra.

\begin{proposition}\label{existence}
Let $A$ be any $F$-algebra and $W$ any $A$-module.
There exists an $F$-linear map
\[ \rho:\vw \rightarrow \End_F(W) \]
with the property that $\rho(\delta)(aw)=a\rho(\delta)(w)+\delta(a)w$
for all $\delta$ in $\vw$, $a$ in $A$ and $w$ in $W$. 
\end{proposition}
\begin{proof} Assume that $g(\delta)=0$. Then there exists a map
$\rho(\delta)$ in $\Hom_F(W,W)$ with the property that $d^0\rho(\delta)=
C(\delta)$. This is if and only if $\rho(\delta)(aw)=a\rho(\delta)(w)+
\delta(a)w$, hence for all $\delta$ in $\vw$ we get a map $\rho(\delta)$,
and the assertion follows.
\end{proof}

Given any $A$-module $W$, we now pick any map $\rho:\vw\rightarrow \End_F(W)$
with the property that $\rho(\delta)(aw)=a\rho(\delta)(w)+\delta(a)w$, which
exists by proposition \ref{existence}. Put $P=\Hom_F(\vw,\End_A(W))$ in sequence
\ref{ex1} and consider the element
$L$ in $\Der_F(A,P)$ defined as follows:
$L(a)(\delta)(w)=a\rho(\delta)(w)-\rho(a\delta)(w)$.
\begin{definition} Let $lc(W)=\overline{L}$ in 
$\HH^1(A,P)=\Ext^1_A(\vw,\End_A(W))$.
\end{definition}
One verifies that the class $lc(W)$ is independent of choice of map
$\rho$, hence it is an invariant of $W$.

\begin{theorem} Let $A$ be any $F$-algebra, and $W$ any $A$-module,
then $lc(W)=0$ if and only if $W$ has a $\vw$-connection.
\end{theorem}
\begin{proof} Assume $lc(W)=0$. Then there exists a map $\eta$ in
$\Hom_F(\vw,\End_A(W))$ with the property that $d^0\eta=L$. Then
the map $\rho+\eta=\nabla:\vw\rightarrow \End_A(W)$ is a $\vw$-connection,
and the assertion follows.
\end{proof} 

From a groupoid in schemes (roughly speaking an algebraic stack, see \cite{LAUMON}) one constructs a  Lie-Rinehart algebra in a way 
similar to the way one constructs
the Lie-algebra from a group-scheme. A natural problem is to find necessary and sufficient criteria for 
the linear  Lie-Rinehart algebra to be integrable to a groupoid in schemes. 

Note that for any $A$-submodule and $k$-sub-Lie algebra $\lg$ of $\Der_k(A)$, there exists a
\emph{generalized universal enveloping algebra} $\U (A,\lg)$ which is a sub-algebra of $\D (A)$, where $\D (A)$ is the ring
of differential-operators of $A$. The algebra $\U (A,\lg)$ has the property that there is a one-to-one correspondence
between $A$-modules with a flat $\lg$-connection and left $\U (A,\lg)$-modules. There exists a generalized PBW-theorem
for the algebra $\U (A,\lg)$ when $\lg$ is a projective $A$-module (see \cite{RINE}). The dual algebra $\U (A,\lg)^*$ is commutative
and its spectrum $\Spec (\U(A,\lg)^*)$ is a formal equivalence-relation in schemes.
Note also that all constructions in this section globalize.

\section{Explicit examples: Algebraic $\vw$-connections}
In this section we apply the theory developed in the previous section
to compute explicitly algebraic $\vw$-connections on a class of maximal
Cohen-Macaulay modules on isolated hypersurface singularities
$B_{mn2}=x^m+y^n+z^2$. Let in the following $F$ be a field of characteristic 
zero, and $A=F[[x,y,z]]/x^m+y^n+z^2$. We are interested in maximal Cohen-Macaulay
modules on $A$, and such modules have a nice description, due to \cite{EIS}:
Consider the two matrices 
\[ \phi=\begin{pmatrix} x^{m-k} & y^{n-l} & 0 & z \\
 y^l & -x^k & z & 0 \\
z & 0 & -y^{n-l} & -x^k \\
0 & z & x^{m-k} & -y^l 
\end{pmatrix} \]
and 
\[ \phi=\begin{pmatrix} x^k & y^{n-l} & z & 0 \\
y^l & -x^{m-k} & 0 & z \\
0 & z & -y^l & x^k \\
z & 0 & -x^{m-k} & -y^{n-l} \end{pmatrix} ,\]
where $1\leq k\leq m$ and $1\leq l \leq n$. Let $f$ be the polynomial $x^m+y^n+z^2$. 
The matrices $\phi$ and $\psi$ have the property that $\phi\psi=\psi\phi=fI$ where
$I$ is the rank 4 identity matrix. Hence we get a complex of $A$-modules
\begin{equation}\label{complex}  \cdots
\rightarrow^\psi P \rightarrow ^\phi P \rightarrow ^\psi P \rightarrow ^\phi P \rightarrow
W(\phi,\psi) \rightarrow 0 .\end{equation}
Note that the sequence \ref{complex} is a complex since $\phi\psi=\psi\phi=fI=0$.
By \cite{EIS}, the module $W=W(\phi,\psi)$ is a maximal Cohen-Macaulay module on $A$.
The ordered pair $(\phi,\psi)$ is a \emph{matrix-factorization} of the polynomial $f$.
We want to compute explicitly algebraic $\vw$-connections on the modules $W=W(\phi,\psi)$
for all $1\leq k\leq m$ and $1\leq l \leq n$, i.e we want to give explicit
formulas for $A$-linear maps $\nabla^{\phi,\psi}=\nabla:\vw\rightarrow \End_F(W)$
satisfying $\nabla(\delta)(aw)=a\nabla(\delta)(w)+\delta(a)w$ for all $a$ in $A$, $w$ in $W$
and $\delta$ in $\vw$. Hence first we have to compute generators and syzygies of 
the derivation-modules $\Der_F(A)$. A straight-forward calculation shows that $\Der_F(A)$ is
generated by the derivations
\[ \delta_0 =2nx\partial_{x} +2my\partial_{y}+mnz\partial_{z} \]
\[ \delta_1 =mx^{m-1}\partial_{y}-ny^{n-1}\partial_{z} \]
\[ \delta_2 =-2z\pa_{x} +mx^{m-1}\pa_{z} \]
\[ \delta_3 =-2z\pa_{y} +ny^{n-1}\pa_{z} ,\]
hence we get a surjective map of $A$-modules 
$\eta:A^{4}\rightarrow \Der_F(A)\rightarrow 0$. A calculation showas that the 
syzygy-matrix of $\Der_F(A)$ is the following matrix
\[
\rho=\begin{pmatrix} y^{n-1} & z & 0 & x^{m-1} \\
 2x & 0 & -2z & -2y \\
0 & nx & ny^{n-1} & -nz \\
-mz & my & -mx^{m-1} & 0 \end{pmatrix} ,\]
hence we get an exact sequence of $A$-modules
\[ \cdots \rightarrow A^4 \rightarrow ^\rho A^4 \rightarrow^\eta \Der_F(A)\rightarrow 0.\]
A calculation shows that the Kodaira-Spencer map $g:\Der_F(A)\rightarrow \Ext^1_A(W,W)$ is
zero for the modules $W$, hence $\vw=\Der_F(A)$, and the calculation also provides us
with elements $\nabla(\delta_i)$ in $\End_F(W)$ with the property that
$\nabla(\delta_i)(aw)=a\nabla(\delta_i)(w)+\delta_i(a)w$ for $i=0,..,3$. Hence we get
an $F$-linear map $\nabla :\vw\rightarrow \End_F(W)$. Given any map $e_i$ in $\End_A(W)$, 
it follows that the map $\nabla(\delta_i)+e_i$ again is an element in $\End_F(W)$ with desired
derivation-property, hence we seek endomorphisms $e_0,..,e_3$ in $\End_A(W)$ with the property
that the adjusted map $\overline{\nabla}:\vw\rightarrow \End_F(W)$ defined by
$\overline{\nabla}(\delta_i)=\nabla(\delta_i)+e_i$ is $A$-linear. We see that we have to solve
equations in the ring $\End_A(W)$. In the examples above, it turns out if one looks closely
that one can find elements in $\End_A(W)$ by inspection so as to reduce the problem to
solve linear equations in the field $F$. If one does this, one arrives at the following
expressions:
\begin{equation} \label{nab0}
\overline{\nabla}_{\delta_0}=\delta_0+A_0= \end{equation}
\[\delta_0+
\begin{pmatrix} nk+ml-\frac{1}{2}mn & 0 & 0 & 0 \\
0 & \frac{3}{2}mn-ml-nk & 0 & \\
 0 & 0 & \frac{1}{2}mn+ml-nk & 0 \\
0 & 0 & 0 & \frac{1}{2}mn+nk-ml \end{pmatrix}.\]

\begin{equation}\label{nab1}
 \onab_{\delta_1}=\delta_1+
\begin{pmatrix} 0 & b_2 & 0 & 0 \\
b_1 & 0 & 0 & 0 \\
0 & 0 & 0 & b_4 \\
0 & 0 & b_3 & 0 \end{pmatrix}=\delta_1+A_1,\end{equation}
with  $b_1=\frac{1}{4}(mn-2nk-2ml)x^{k-1}y^{l-1}$, $b_2=\frac{1}{4}(3mn-2ml-2nk)
x^{m-k-1}y^{n-l-1}$, $b_3=\frac{1}{4}(2nk-mn-2ml)x^{m-k-1}y^{l-1}$ and
$b_4=\frac{1}{4}(2nk-2ml+mn)x^{k-1}y^{n-l-1}$.

\begin{equation}\label{nab2}
 \onab_{\delta_2}=\delta_2+
\begin{pmatrix} 0 & 0 & c_3 & 0 \\
0 & 0 & 0 & c_4 \\
c_1 & 0 & 0 & 0 \\
0 & c_2 & 0 & 0 \end{pmatrix}=\delta_2+A_2 ,\end{equation}
with $c_1=\frac{1}{n}(\frac{1}{2}mn-ml-nk)x^{k-1},$ $c_2=\frac{1}{n}(
\frac{3}{2}mn-ml-nk)x^{m-k-1}$, $c_3=\frac{1}{n}(\frac{1}{2}mn+ml-nk)x^{m-k-1}$
and $c_4=\frac{1}{n}(ml-nk-\frac{1}{2}mn)x^{k-1}$.

\begin{equation}\label{nab3}
 \onab_{\delta_3}=\delta_3+
\begin{pmatrix} 0 & 0 & 0 & d_4 \\
0 & 0 & d_3 & 0 \\
0 & d_2 & 0 & 0 \\
d_1 & 0 & 0 & 0 \end{pmatrix}=\delta_3+A_3,\end{equation}

where $d_1=\frac{1}{m}(\frac{1}{2}mn-ml-nk)y^{l-1}$, $d_2=\frac{1}{m}(ml+nk
-\frac{3}{2}mn)y^{n-l-1}$, $ d_3=\frac{1}{m}(\frac{1}{2}mn+ml-nk)y^{l-1}$ and 
$d_4=\frac{1}{m}(\frac{1}{2}mn-ml+nk)y^{n-l-1}$.

\begin{theorem} \label{brieskorn}
For all $1\leq k \leq m$ and $1\leq l \leq n$ the equations
\ref{nab0}-\ref{nab3} define algebraic
$\vw$-connetions $\overline{\nabla}^{\phi,\psi}:\vw\rightarrow \End_F(W)$ where $W=W(\phi,\psi)$.
\end{theorem}
\begin{proof} The module $W$ is given by the exact sequence
\[ \cdots \rightarrow ^\psi A^4 \rightarrow ^\phi A^4 \rightarrow W \rightarrow 0 ,\]
hence an element $w$ in $W$ is an equivalence class $\overline{a}$ of an element
$a$ in $A^4$. Let $a$ in $A^4$ be the element
\[ a= \begin{pmatrix}  a_1 \\ a_2 \\ a_3 \\ a_4 \end{pmatrix} \]
and consider the class $w=\overline{a}$ in $W=A^4/im\phi$. Define the $\vw$-connection
$\overline{\nabla}$ as follows: 
\[ \overline{\nabla}(\delta_i)(w)=(\delta_i+A_i)
\begin{pmatrix} a_1 \\ a_2 \\ a_3 \\ a_4 \end{pmatrix} .\]
Then one verifies that this definition gives a well-defined $A$-linear map
$\overline{\nabla}^{\phi,\psi}=\overline{\nabla}:\vw\rightarrow \End_F(W)$, and we 
have proved the assertion.
\end{proof}
Note that the connections $\overline{\nabla}^{\phi,\psi}$ from Theorem \ref{brieskorn}
exist over any field $F$ of characteristic prime to $m$ and $n$.

In Theorem \ref{brieskorn} we saw examples where $\vw$ was the whole module of derivations for a class of modules on
some hypersurface singularities. In \cite{maa10},
Theorem 5.1 the splitting type of the principal parts $\Pr^k(\O(d))$ is calculated on the projective line over a field
of characteristic zero. When $d\geq 1$, he following formula is shown:
\[ \Pr^1(\O(d))\cong \O(d-1)\oplus \O(d-1) . \]
It follows that the Atiyah sequence (see \cite{ATIYAH}) 
\[ 0\rightarrow \Omega^1\otimes \O(d)\rightarrow \Pr^1(\O(d)) \rightarrow \O(d) \rightarrow 0 \]
does not split hence $\O(d)$ does not have a connection 
\[ \nabla :\O(d) \rightarrow \Omega^1\otimes \O(d) .\]
It follows that there does not exist an action
\[ \rho: T_{\P^1}\rightarrow \End(\O(d)) , \]
hence we see that for $\O(d)$ on $\P^1$ over a field of characteristic zero, 
the linear  Lie-Rinehart algebra  $\mathbf{V}_{\O(d)}$ is a strict sub-sheaf
of the tangent sheaf $T_{\P^1}$. 

\section{Chern-classes}
In this section we define for any $F$-algebra $A$ where $F$ is any field, 
and any $A$-module $W$
of finite presentation with a $\vw$-connection the first Chern-class
$c_1(W)$ in $\H^2(\vw|_U,\O_U)$, where $U$ in $\Spec(A)$ is the open subset
where $W$ is locally free, and $\H^i(\vw|_U,\O_U)$ is the Chevalley-Hochschild
cohomology of the restricted  Lie-Rinehart algebra $\vw|_U$ with values in the sheaf
$\O_U$. 

\begin{definition} Let $\lg$ be a  Lie-Rinehart algebra and $W$ a $\lg$-module.
The Chevalley-Hochshild complex $\C^\bullet(\lg,W)$ is defined as follows:
\[ \C^p(\lg,W)=\Hom_A(\lg^{\wedge p},W) ,\]
with differentials $d^p:\C^p(\lg,W)\rightarrow \C^{p+1}(\lg,W)$ defined by
\[ d^p\phi(g_1\wedge \cdots \wedge g_{p+1})=\sum_{i=1}^{p+1}(-1)^{i+1}
g_i\phi(g_1\wedge \cdots \wedge \overline{g_i}\wedge \cdots \wedge g_{p+1})+\]
\[\sum_{i<j}(-1)^{i+j}\phi([g_i,g_j]\wedge \cdots \wedge\overline{g_i}\wedge \cdots
\wedge \overline{g_j}\wedge \cdots \wedge g_{p+1}).\]
Here $g\phi(g_1\wedge \cdots \wedge g_p)=\nabla(g)\phi(g_1\wedge \cdots \wedge g_p)$
and overlined elements should be excluded. The cohomology $\H^i(\C^\bullet(\lg,W))$
is denoted $\H^i(\lg,W)$.
\end{definition}
Note that $\C^\bullet(\lg,W)$ is a complex if and only if the $\lg$-connection
$\nabla$ is \emph{flat}, i.e a morphism of Lie-algebras. The complex $\C^\bullet(\lg,W)$
is a complex generalizing simultaneously the algebraic deRham-complex and Chevalley-Eilenberg complex.
Now given a $\vw$-connection $\nabla$ on $W$, we immediately get a $\vw$-connection
on $\End_A(W)$, denoted $ad\nabla$. We see that if $\nabla$ is flat, it
follows that $ad\nabla$ is flat. Consider the open subset $U$ of $\Spec(A)$
where $W$ is locally free. We restrict the $\vw$-connection to $U$, to get
a connection $\nabla:\vw|_U\rightarrow \End_{\O_U}(W|_U)$. Since $W|_U$ is locally 
free we can find an open cover ${U_i}$ of $U$ where $U_i=D(f_i)$ and 
$W|_{D(f_i)}=A_{f_i}^n$. On the open subsets $U_i$ we get restricted connections
$\nabla_i$, but since $W|_{U_i}$ is a free module, we have on $U_i$ when we pick a
basis a flat connection $\rho_i$. Therefore on any open subset $U_i$ we can 
consider the complex of sheaves $\C^\bullet(\vw|_{U_i},\End_{\O_{U_i}}(W|_{U_i}))$ with respect
to the regular connection $ad\rho_i$. We see that $R_\nabla|_{U_i}$ 
is an element of $\C^2(\vw|_{U_i},\O_{U_i})$. One checks that the element
$c^i=trace \circ R_\nabla|_{U_i}$ is a cocycle of the complex 
$\C^\bullet(\vw|_{U_i},\O_{U_i})$ for all $i$. We have constructed elements
\[ c^i \in \C^2(\vw|_{U_i},\O_{U_i}) \]
which coincide on intersections since trace is independent with respect to choice of 
basis, hence the sheaf-structure on $\C^2(\vw|_{U},\O_{U})$ gives a uniquely
defined element $c$ in $\C^2(\vw|_{U},\O_{U})$, such that $c|_{U_i}=c^i$
for all $i$. There exists a commutative diagram
\[ \diagram \C^p(\vw|_U,\O_U) \rto^{d^p} \dto^{|_{U_i}} & \C^{p+1}(\vw|_U,\O_U)
\dto^{|_{U_i}} \\
\C^p(\vw|_{U_i},\O_{U_i}) \rto^{d^p} & \C^{p+1}(\vw|_{U_i},\O_{U_i})
\enddiagram ,\]
which proves that the element $c$ is a cocycle in the complex $\C^\bullet(\vw|_U,
\O_U)$. 

\begin{theorem} \label{chernclass}
There exists a class $c_1(W)$ in $\H^2(\vw|_U,\O_U)$ which 
is independent with respect to choice of $\vw$-connection.
\end{theorem}
\begin{proof} Existence of the class $c_1(W)$ follows from the argument above.
Independence with respect to choice of connection is straightforward.
\end{proof}
Note that if the $\vw$-connection is flat, the first Chern-class $c_1(W)$ is zero.
Note also that the construction in this section cam be done with any $A$-module $W$ of finite presentation with
a $\lg$-connection, where $\lg$ is any  Lie-Rinehart algebra. 
\begin{theorem} \label{zero}
The $\vw$-connections $\overline{\nabla}^{\phi,\psi}$ calculated in Theorem \ref{brieskorn}
 are flat hence $c_1(W(\phi,\psi))=0$ for all the modules $W(\phi,\psi)$ on the singularities
$B_{mn2}$.
\end{theorem}
\begin{proof} Easy calculation.
\end{proof}
Note that the flat $\vw$-connections in Theorem \ref{zero} give rise to a class of left modules on the algebra
of differential operators $\D (A)$ where $A$ is the ring $k[x,y,z]/f$ and  $f=x^m+y^n+z^2$. 
Note also that Kohno has in \cite{KOHNO} computed 
the Alexander-polynomial of an irreducible plane curve $C$ in $\mathbf{C}^2$ using a certain
logarithmic deRham-complex $\Omega^\bullet_{\mathbf{C}^2}(*C)$. It would be interesting to see if 
the Alexander-polynomial can be computed in terms of a $\lg$-connection.

\section{Surface quotient-singularities}
In this section we consider maximal Cohen-Macaulay modules on quotient singularities 
of dimension two, and their first Chern-class. Let now $F=\mathbf{C}$ be the
complex numbers, and let $G\subseteq GL(2,F)$ be a finite sub-group with no
pseudo-reflections. Consider the natural action $G\times F^2\rightarrow F^2$, and
the quotient $X=F^2/G$. It is an affine scheme with an isolated singularity. Pick
a representation $\rho:G\rightarrow GL(V)$ where $V$ is an $F$-vectorspace,
and consider the $A$-module $V\otimes_FA$ where $A=F[x,y]$. Define
a $G$-action as follows: $g(v\otimes a)=\rho(g)v\otimes ga$, and let
$W_\rho=(V\otimes_F V)^G$ be the $G$-invariants of the $G$-action defined.
Then by \cite{HER} $W_\rho$ is a maximal Cohen-Macaulay module on $A^G$. If $\rho$
is indecomposable, $M_\rho$ is irreducible. The McKay correspondence in dimension 2
says that all maximal Cohen-Macaulay modules
on $A^G$ arise this way. Define a $G$-action on $\lg=\Der_F(A)$ as follows:
$(g\delta)(a)=g\delta(g^{-1}a)$, and let $\lg^G$ be the $G$-invariant
derivations. The module $V\otimes_F A$ is a free $A$-module, hence there exists
trivially a regular $\lg$-connection on $V\otimes_F A$. This implies that
we get an induced $\lg^G$-connection on $V\otimes_F A$.

\begin{proposition} There exists an action of $\lg^G$ on $W_\rho$.
\end{proposition}
\begin{proof} Straightforward.
\end{proof}

From \cite{SCH} it follows that $\lg^G$ is isomorphic to $\Der_F(A^G)$, hence
we have proved that all maximal Cohen-Macaulay modules $W_\rho$ on 
$\Spec(A^G)$ posess a $\Der_F(A^G)$ -connection, and these are all regular
connections since they are induced by the trivial one on $V\otimes_F A$.

\begin{theorem} \label{quotient}
Let $X=\Spec(A^G)$ be a 2-dimensional quotient singularity
and let $W_\rho$ be a maximal Cohen-Macaulay module on $X$ then $c_1(W_\rho)=0$.
\end{theorem}
\begin{proof} Follows from the argument above.
\end{proof}

\textbf{Acknowledgements} This paper is a slight extension of my  master-thesis written 
at the University of Oslo under supervision of A.O Laudal, and I want to thank him 
for suggesting the problems discussed in this paper. Thanks also to 
J. Cristophersen for pointing out to me the existence of the connections in section 5.

\end{document}